\documentclass[fleqn, leqno, 12pt]{article}
\usepackage{amsmath, amssymb, bm}

\newcommand{\thmCRhopf}{Theorem 2.1}
\newcommand{\propmain}{Proposition 2.1}
\newcommand{\thmmain}{Theorem 2.2}
\newcommand{\remcomm}{Remark 2.1}
\newcommand{\remcompa}{Remark 2.2}
\newcommand{\remmer}{Remark 2.3}

\newcommand{\thmbasic}{Theorem 3.1}
\newcommand{\remext}{Remark 3.1}
\newcommand{\remconc}{Remark 3.2}

\newcommand{\propisom}{Proposition 3.1}
\newcommand{\propkernel}{Proposition 3.2}
\newcommand{\proppsd}{Proposition 3.3}

\newcommand{\lemsi}{Lemma 4.1}
\newcommand{\lemti}{Lemma 4.2}

\newcommand{\lemfol}{Lemma 5.1}
\newcommand{\lembas}{Lemma 5.2}

\newcommand{\Hone}{H1}
\newcommand{\Htwo}{H2}
\newcommand{\Hthree}{H3}
\newcommand{\Hfour}{H4}
\newcommand{\Bone}{B1}
\newcommand{\Btwo}{B2}
\newcommand{\Bthree}{B3}
\newcommand{\Kone}{K1}
\newcommand{\Ktwo}{K2}
\newcommand{\KtwoFirst}{K2-1}
\newcommand{\KtwoSecond}{K2-2}
\newcommand{\KtwoThird}{K2-3}
\newcommand{\KtwoFourth}{K2-4}

\newcommand{\eqpde}{1.1}

\newcommand{\eqpersp}{2.1}
\newcommand{\eqpu}{2.2}
\newcommand{\eqpuu}{2.3}
\newcommand{\eqtpe}{2.4}

\newcommand{\eqH}{3.1}
\newcommand{\eqThe}{3.2}
\newcommand{\eqgi}{3.3}
\newcommand{\eqlon}{3.4}
\newcommand{\eqtoa}{3.5}

\newcommand{\eqfrac}{4.1}
\newcommand{\eqxn}{4.2}
\newcommand{\eqyn}{4.3}
\newcommand{\eqmapdh}{4.4}
\newcommand{\eqgam}{4.5}
\newcommand{\eqmuz}{4.6}
\newcommand{\eqfly}{4.7}
\newcommand{\eqsvec}{4.8}
\newcommand{\eqrto}{4.9}
\newcommand{\eqluo}{4.10}
\newcommand{\eqpin}{4.11}
\newcommand{\eqwg}{4.12}
\newcommand{\eqwt}{4.13}
\newcommand{\eqndn}{4.14}
\newcommand{\eqnn}{4.15}
\newcommand{\eqdnu}{4.16}

\newcommand{\eqhc}{5.1}
\newcommand{\eqhlp}{5.2}
\newcommand{\eqsp}{5.3}
\newcommand{\equni}{5.4}
\newcommand{\eqint}{5.5}
\newcommand{\eqfgh}{5.6}
\newcommand{\eqplp}{5.7}
\newcommand{\eqazi}{5.8}
\newcommand{\equh}{5.9}
\newcommand{\equhz}{5.10}
\newcommand{\eqazz}{5.11}
\newcommand{\eqzin}{5.12}
\newcommand{\eqfns}{5.13}
\newcommand{\eqini}{5.14}
\newcommand{\eqiac}{5.15}
\newcommand{\equva}{5.16}
\newcommand{\eqvmp}{5.17}
\newcommand{\eqvne}{5.18}
\newcommand{\eqnux}{5.19}
\newcommand{\eqnov}{5.20}

\newcommand{\refAma}{A}
\newcommand{\refBKST}{BKST}
\newcommand{\refCRbif}{CR1}
\newcommand{\refCRhopf}{CR2}
\newcommand{\refGMW}{GMW}
\newcommand{\refKsym}{K1}
\newcommand{\refKcom}{K2}
\newcommand{\refKsome}{K3}
\newcommand{\refKie}{Kie}
\newcommand{\refLiZY}{LiZY}
\newcommand{\refliu}{LMR}
\newcommand{\refMS}{MS}
\newcommand{\refNishi}{NTY}
\newcommand{\refWYZ}{WYZ}

\newcommand{\ubold}{{\bm u}}
\newcommand{\usbold}{{\bm u}_\star}
\newcommand{\etabold}{{\bm \eta}}

\newcommand{\psis}{\psi_\star}
\newcommand{\psibold}{{\bm \psi}}
\newcommand{\psisbold}{{\bm \psi}_\star}

\newcommand{\cV}{ {\cal V} }
\newcommand{\cX}{ {\cal X} }
\newcommand{\cY}{ {\cal Y} }

\newcommand{\sech}{\text{sech}\,}

\newcommand{\propspace}{{\vskip 0.4 cm}\n}
\newcommand{\eqspace}{{\vskip 0.25 cm}\n}
\def\tag#1 {\leqno{(#1)} }

\newcommand{\Vc}{V_{\hbox{c}}}
\newcommand{\Uc}{U_{\hbox{c}}}
\newcommand{\Ac}{A_{\hbox{c}}}

\newcommand{\us}{u_\star}
\newcommand{\eone}{{\bm e}_1}

\DeclareMathOperator{\Realpart}{Re}
\DeclareMathOperator{\Imaginarypart}{Im}

\newcommand{\Azc}{A_{0 {\hskip 0.03cm} c}}
\newcommand{\Ai}{A_{\i}}
\newcommand{\Aic}{A_{{\i} {\hskip 0.03cm} c}}

\def\linear{ {\cal L} }
\def\domain{ {\cal D} }

\def\range{ {\cal R} }
\def\kernel{ {\cal N} }

\def\span{\hbox{span}}

\newcommand{\directsum}{\oplus}

\newcommand{\codim}{\text{codim}\,}

\def\q{\quad}

\def\text{\hbox}
\def\medn{\medbreak\noindent}
\def\n {\par\noindent}

\def\newchapter{\vskip 1.0cm}

\def\END{ \hfill \hbox{ \vrule
     \vbox{\hrule  \vskip 3pt
          \hbox{\hskip 3pt
               \vbox{\hsize = .7in \raggedright
                          }%
                   \hskip 3pt}%
                \vskip 3pt \hrule}%
             \vrule}
 }

\def\moji#1 {{\hskip 0.25cm}\hbox{#1}{\hskip 0.25cm}}

\def\; {\,;\,}
\def\: {\,:\,}

\def\N{ \mathbb{N} }
\def\R{ \mathbb{R} }
\def\Z{ \mathbb{Z} }
\def\C {\mathbb{C}}

\def\e{\varepsilon}

\def\a{\alpha}
\def\b{\beta}
\def\lam{\lambda}
\def\Lam{\Lambda}
\def\i{\infty}

\begin{document}


\centerline{ \bf  The Hopf bifurcation theorem in Hilbert spaces}
\centerline{ \bf  for abstract  semilinear equations}

  \vskip 0.3 cm

 \centerline{Tadashi KAWANAGO}
\vskip 0.3 cm
\centerline{Faculty of Education, Saga University, Saga 840-8502, Japan}
\centerline{e-mail: tadashi@cc.saga-u.ac.jp}

 \vskip 1.5 cm

{\bf Abstract}
\n		
We prove a Hopf bifurcation theorem in Hilbert spaces
for abstract semilinear equations, which  
 improves a  classical result by 
Crandall and Rabinowitz in the case where basic spaces are 
Hilbert spaces. Actually, our theorem does not need any compactness
conditions, which leads to wider applications. 
In particular, our  theorem can be applied to
semilinear equations in unbounded domains of $\R^n$.

\newchapter

\centerline{\bf 1. Introduction}
\medn
Concerning the Hopf bifurcation theorems in infinite dimensions,
a lot of versions have been proved until now
 (see e.g.   [\refCRhopf],  [\refAma],  [\refliu],
[\refGMW] and the references therein).
Among them [\refCRhopf, Theorem 1.11] by Crandall and Rabinowitz is 
one of most important results. It is a theorem for
abstract  semilinear equations
and has been well applied so far
to various studies
because of its generality (see e.g. [\refGMW] and [\refWYZ]). 
 It needs, however,
  some compactness condition, and, consequently, can not
  be applied   to partial differential equations
  in unbouded domains  of $\R^n$.

On the other hand, Hopf bifurcation in partial differential equations 
  in the unbounded domain of $\R^n$ has been  studied 
  more recently and
   Hopf bifurcation theorems applicable to
 such studies were proven
(see e.g. [\refLiZY], [\refMS] and [\refBKST]). 
As far as the author knows, however, 
each of them can be
applied  to a specific type of eqautions, to be sure,
but it does not have generality  applicable to various studies.  

          In this paper we prove a general Hopf bifurcation 
theorem  for abstract  nonlinear  \linebreak
equations 
which is applicable to semilinear differential eqautions
 in unbouded domains of $\R^n$.
 Actually, our theorem can be considered to be an imporvement of        
          [\refCRhopf, Theorem 1.11] 
in the case where basic spaces are Hilbert spaces. 

			We consider the next abstract semilinear equation 
in Hilbert spaces:
\[
		u_t = A u + h(\lam, u),			\tag{\eqpde}
\]
where the linear operator $A$ and the map  $h$ are described 
in Section 2 below. 

		The assumptions of our main theorem ({\thmmain} below) are
weaker than those of   [\refCRhopf, Theorem 1.11].
Actually, our result has the following features:
 \begin{itemize}
  \item 
We do not assume  that  $A$ generates a $C_0$-semigroup.
  \item
We do not assume  that 
 $A$ has compact resolvents.
   \end{itemize}
These features  contribute to  wider applications (see Section 5 below).
More concretely, the former  has merit since Hopf bifurcation can 
occur even for the case where $A$ is not densely defined
 (see [\refliu] and the references therein), for example. The latter feature makes it possible to apply
our {\thmmain}  to semilinear equations on unbounded domains of $\R^n$.
Actually, we treat the Cauchy problem for  a system of semilinear heat equations as a concrete example in Section 5 below.

			The plan of our paper is the following.
In Section 2 we describe our main results and discuss the features of our
results. In Section 3 we describe some preliminary results to prove our main results. We prove our main result in Section 4. In Section 5
we present some concrete examples. 
 
 \newchapter

\centerline{\bf 2. The Hopf bifurcation theorems in infinite dimensions}

			In this section we present a new bifurcation theorem
({\thmmain} below), which is an infinite dimensional version of 
the classical Hopf bifurcation theorem.
 
		Let $V$ be a real Banach space and $V_c= V + i V$ be  its 
complexification.
Let $A$ be a closed linear operator on $V$ 
with a bounded inverse $A^{-1}$.
We denote its domain by $\domain(A)$, range by $\range(A)$,
null space by $\kernel(A)$ and the extension of $A$ on $V_c$ 
by $\Ac$. If $W$ is another Banach space,
$\linear(V, W)$ denotes the set of bounded linear operators from
$V$ to $W$. We simply write $\linear(V):= \linear(V, V)$.

			We consider the  equation (\eqpde).
First we describe a known result 
{  [\refCRhopf, Theorem 1.11]}.
We assume the following conditions (\Hone) - (\Hfour).

\propspace
{(\Hone)}  The operator $A$ is the generator of a 
$C_0$-semigroup on $V$.
\propspace
{(\Htwo)}   $\exp(t \Ac)$ is a holomorphic semigroup  on $V_c$.
\propspace
{(\Hthree)}  
The resolvent $(z - A_c)^{-1}$ is {  compact} 
for any $z \in \rho(A_c)$.

\propspace
It follows from (\Hone) and (\Htwo) that if $r > \Realpart z$ 
for all $z \in \sigma(A_c)$, then the fractional powers $(r - A)^{\a}$
are defined for  $\a \ge 0$. 
We can define the Banach space 
$V_{\a} \subset V$
with norms $\| \cdot \|_{\a}$  by 
$V_{\a} := \domain((r - A)^{\a}) \moji{with} 
\| v \|_{\a} := \| (r - A)^{\a} v \|_{V} \moji{for} v \in V_{\a}$.

\propspace
{(\Hfour) } There exist an $\a \in [0, 1)$ and an open 
neighborhood $\Omega$ of $(0,0)$
 in $\R \times V_{\alpha}$
such that $h \in C^2(\Omega, V)$. Moreover, $h_u(0, 0)=0$ and
$h(\lam, 0) = 0$   if  $(\lam, 0) \in \Omega$.
{\vskip 0.5cm}
We also assume the following (\Bone) - (\Bthree):

\propspace
(\Bone)  $\,\,\,\pm\, i$ are the  simple eigenvalues of  $A_c$, i.e. 
\[
\begin{cases}
& \dim \kernel (i - A_c) = 1  =\codim \range(i - A_c), \hfill  \\ 
& \psi \in \kernel(i - A_c) - \{ 0 \} \,\,
	\Longrightarrow \,\, \psi \not \in \range(i-A_c). \hfill \\ 
\end{cases}
\]
So, by the implicit function theorem,
  $A_c + h_u(\lam, 0)$ has an eigenvalue $\mu(\lam) \in \C$
and eigenfunction $\psi(\lam) \in \domain(A_c)$
corresponding to $\mu(\lam)$  
for any $\lam$ in a small  neighborhood of $0$ 
such that $\mu(0)= i$ and that 
$\mu(\lam)$ and $\psi(\lam)$
 are functions  of class $C^2$.
 
\propspace
(\Btwo)\q	(Transversality condition of eigenvalues)\q 
$\text{Re}\,\mu'(0) \not= 0$.

\propspace
(\Bthree) \q $i k \in \rho(A_c)$
		for  $k\in \Z - \{ -1, 1 \}$.

\propspace
{\bf  \thmCRhopf}. ({  [\refCRhopf, Theorem 1.11]})
 {\sl  We assume  (\Hone) - (\Hfour)  and (\Bone) - (\Bthree). 
 Then, $(\lam, u)= (0,0)$ is a Hopf bifurcation point of (\eqpde).
}
\propspace
\n
Here, we omit the description in [\refCRhopf, Theorem 1.11]
on the uniqueness of the branch of bifurcating periodic solutions.

			Next, we state our new results.
We consider the case in which 
$V$ is a real Hilbert space and $0 \in \rho(A_c)$.
We define the real Hilbert space $U:=\domain(A) \subset V$
with the norm $\| u \|_U := \|  A u \|_{V}$ for $u \in U$.

		We set the  real Hilbert spaces $X$ and $Y$ by
\[
  X:= H^1_{\text{per}} ((0, 2\pi), V) \cap L^2((0, 2\pi), U)
 \moji{and}  Y:= L^2((0, 2\pi), V).			\tag{\eqpersp}
\]
Here, $H^1_{\text{per}} ((0, 2\pi), V):= 
	\{ u \in H^1 ((0, 2\pi), V) \; u(0)=u(2\pi) \}$.

		We assume (\Bone) - (\Bthree) and the following 
(\Kone), (\KtwoFirst) - (\KtwoFourth):

\propspace
\n
 (\Kone) \hspace{0.4 cm}  There exists $M \in (0, \infty)$ such that
\[
\| (i n - A_c)^{-1} \|_{V_c \to V_c} 
		\le \dfrac{M}{n} \moji{for} n= 2, 3, 4,\cdots.
\]
(\KtwoFirst) {\hskip 0.2cm} 
 There is an open interval $K$ in $\R$
 such that $0 \in K$ and $h$ is a map from $K \times U$ to $V$.
\propspace
For any $(\lam, u) \in K \times X$, we set 
$[h(\lam, u)](t) := h(\lam, u(t))$ for a.e. $t \in (0, 2\pi)$.
 \propspace
 \n
(\KtwoSecond) {\hskip 0.2cm} 
$h(\lam, u) \in Y$ for any $(\lam, u) \in K \times X$.

\propspace
We define the map $\Phi \: (\lam, u) \in K \times X 
\mapsto h(\lam, u) \in Y$.
 \[
\Phi \in C^2(K \times X, Y).		\tag{\KtwoThird}
\]

\propspace
{\bf \remcomm}.
We can regard $U$ (resp.\,\,$V$) as the closed subspace of 
$X$ (resp. resp.\,\,$Y)$ which consists of constant functions in $X$
 (resp.\,\,$Y)$.
Then, we verify that (\KtwoThird) implies 
$h \in C^2(K \times U, V)$
with
\[
[\Phi_u (\lam, u) v](t)= h_u (\lam, u(t)) v(t)   \moji{in} V,	\tag{\eqpu}
\]
\[
[\Phi_{uu} (\lam, u) v w](t)= h_{uu} (\lam, u(t)) v(t) w(t)   \moji{in} V
													\tag{\eqpuu}
\]
and so on for $\lam \in K$, $u, v, w\in X$ and a.e. 
$t \in (0, 2\pi)$.
		\END

\propspace
(\KtwoFourth) {\hskip 0.2cm} 
$h_u(0, 0)=0$ and $h(\lam, 0) = 0$   if  $\lam \in K$.
{\vskip 0.25cm}
\n
In what follows we simply denote (\KtwoFirst) - (\KtwoFourth) 
by (\Ktwo). Now, we shortly state our result:
\eqspace
{\bf  \propmain}.  
{\sl  Let $V$ be a real Hilbert space and $A$ be a closed linear operator
on $V$.
We assume (\Kone), (\Ktwo) and  (\Bone) - (\Bthree). 
Then, $(\lam, u)= (0,0)$ is a Hopf bifurcation point of (\eqpde).
}

\propspace
{\bf \remcompa}.
The conditions (\Hone) and (\Htwo) imply (\Kone).
Though we do not assume (\Hone) in {\propmain}, we note that
(\Kone) and (\Hone) imply (\Htwo). 
	\END
\propspace

\n
{\propmain} is a short  version of our main result {\thmmain} below, 
which shows that the branch of bifurcating periodic solutions  are
unique in a neighborhood of $(\lam, u)= (0,0)$.

			Next, we make preparation to state our main result.			
Let $m \in \Z$, $n \in \N$  and $u \in V_c$.
We write $e_m(t) := e^{imt}$, $c_n(t):= \cos nt$ and
$s_n(t):= \sin nt$  for $t \in \R$.
We denote $(u \otimes e_m)(t) := u e_m(t) = u e^{imt}$ ($t \in \R$).
 Similarly, $(u \otimes c_n)(t) := u \cos nt$
and $(u \otimes s_n)(t) := u \sin nt$ ($t \in \R$).		
We set
$X_1 := \{ u \otimes c_1 + v \otimes s_1 \; u, v \in U \}$
   as a subspace of $X$. 
We define the translation operator $\tau_\theta$ by
$(\tau_\theta u)(t):= u(t - \theta)$ for any $\,\theta \in \R$.
			For simplicity, we set $f(\lam, u)= Au + h(\lam, u)$.
If $u(t)$ is a $2\pi$-periodic solution of the next equation:
\[
u_t = (\sigma + 1) f(\lam, u)		\tag{\eqtpe}
\]
 then $u(t/(\sigma + 1))$ is a $2\pi (\sigma + 1)$-periodic 
solution of (\eqpde).
   
	 	Our main theorem is the following:

\propspace
{\bf  \thmmain}.  
{\sl We assume all conditions in {\propmain}.
Then, there exist $a, \e > 0$, $u_{\star} \in X_1 - \{ 0 \}$ and
 functions $\zeta =(\lam, \sigma) \in C^1([0, a), \R^2)$, 
$\eta \in C^1([0, a), X)$ with the following properties:
\medn
(a)  \,\, $(\lam, \sigma, u)= (\zeta(\a), \a u_{\star} + \a \eta(\a) )$
is a solution of (\eqtpe) for any $\a \in [0, a)$,
\medn
(b) \,\, $\zeta(0)=\zeta'(0)= (0,0)$ and $\eta(0)=0$,
\medn
(c) \,\, If $(\lam, v)$ is a solution of (\eqpde) of period 
$2\pi(\sigma + 1)$, $|\lam| < \e$, $|\sigma| < \e$,
${\tilde v} \in X$ and
$\|  {\tilde v} \|_X < \e$,  where
${\tilde v}(t):= v((\sigma + 1)t )$ for $t \in \R$, then there exist
$\a \in (0, a)$ and $\theta \in [0, 2\pi)$ such that 
$(\lam, \sigma)= \zeta(\a)$ and $v((\sigma + 1)t)= \a u_{\star}
(t + \theta) + \a \eta(\a) (t + \theta)$ for any $t \in \R$.
}
\propspace
{\bf \remmer}.  (i)
  In our {\thmmain} we do not need the compact resolvent
   condition (\Hthree).
This contributes to wider applications.	
See concrete examples in Section 5.	
\par
		(ii)  We can naturally extend the domain of 
 $\zeta$ and $\eta$ in  {\thmmain} for $\a \in (- a, 0)$
as $C^1$ maps satisfying (a).
 See the proof of {\thmmain} for details.		\END

\propspace
          
           At the end of this section we  explain the feature of our proof
  of {\thmmain}, which is described in Section 4.
The proof of our main result {\thmmain} is based on
 {\thmbasic} below, which is an abstract bifurcation theorem.
           The technique of our proof  differs
from the standard ones used in the proofs of the related known results.
In general the derivation of infinite-dimensional bifurcation theorem
is so far based mainly on two  techniques. 
One is to analyze directly the 
infinite-dimensional space.  
Crandall and Rabinowitz [\refCRhopf]
adapted this way by using the semigroup theory.
Another is to reduce the problem to a finite-dimensional problem
by using the Lyapunov-Schmidt method, the center manifold theorem
and so on (see e.g. [\refKie], [\refliu]).
 Alternatively, we reduce our problem to the  analysis on 
 two infinite-dimensional spaces by using our {\thmbasic} 
 mentioned above. 
 Actually we can express $X$
  (the real space of $2 \pi$-periodic $V$-valued functions defined
   in (\eqpersp)) 
  as the direct sum of a low-frequency subspace
 and a high-frequency subspace. 
 We reduce our problem to the analysis on each subspace.
  Analysis on the low-frequency subspace is 
  as follows:
 The complex space $V_c$ and the real space of $V$-valued 
 simple harmonic oscillation are isomorphic as real linear spaces
 (see {\propisom} below).
 The seemingly difficult points of analysis for the low-frequency subspace
 can be reduced to the linear algebraic properties of the isomorphism.
 This analysis is so simple that  it always works well without
 the choice of functional spaces. Indeed it works well even if
 $V$ is a general Banach space. 
 On the other hand,  the analysis on the high-frequency subspace
  is based on
 the Fourier analysis. Whether it works well or not seems to depend
  much on the choice of functional spaces. We need  here 
  our technical assumption that $V$ is a Hilbert space. 
 
\newchapter

 \centerline{\bf 3. Preliminary results}
\medbreak 

\medbreak
	To begin with, we describe [\refKsome, Theorem 3]
for the case $m=2$.  The proof of our main theorem ({\thmmain})
is based on this result. 

			Let $\cX$ and  $\cY$ be real Banach spaces and 
$O$ be an open neighborhood of $0$ in $\cX$.
Let $J$ be an open neighborhood of $(0, 0)$ in $\R^2$.
Let  $g \in C^2( J \times O, \,\cY)$ be a 
map such that
\[
g(\Lam, 0)= 0 \moji{for \,\, any} 
	\Lam= (\Lam_1, \Lam_2) \in J.
\]
We define $H \: J \times \cX \to \R^2 \times \cY$ by
\[
H
\begin{pmatrix}
\Lam  \\
u  \\
\end{pmatrix}
:=
\begin{pmatrix}
l u - \eone  \\
g_u(\Lam, 0)u  \\
\end{pmatrix}.			\tag{\eqH}
\]
Here, $l = (l^1, l^2) \in \linear(\cX, \R^2)$
and $\eone :=  (1, 0)$.
We define 
$G \: J \times O\, \to \,
		\R\times \cY$ by
\[
G
\begin{pmatrix}
\Lam  \\
u  \\
\end{pmatrix}
:=
\begin{pmatrix}
l ^2 u  \\
g(\Lam, u)  \\
\end{pmatrix}.
\] 
We set  $Z:=  \{ u \in \cX \; l^1 u= 0 \}$. 
\eqspace
{\bf \remext}.
The system such as $H(\Lam, u) = 0$ and $G(\Lam, u) = 0$ above 
is called a {\sl extended system}
in general in the field of numerical analysis 
 (see e.g. [\refKsome]).			\END

\propspace
\n
{\bf \thmbasic}. 
 ([\refKsome, Theorem 3] for the case $m=2$) 
 \n
{\sl In addition to the 
assumptions above we assume that there exists  $\us \in O$ such that
\medn
{(\eqThe)}\q  the extended system $H(\Lam, u) = 0$
has an isolated solution $(\Lam, u)= (0,\us)$. 
\medn
 Then there exist an open neighborhood $W$ of (0,0)  
 in $\R^2 \times \cX$,
$a  \in (0, \infty)$ and  functions 
$\zeta \in C^1((-a, a), \R^2)$, 
$\eta \in C^1((-a, a), Z)$ such that $\zeta(0)=0$, $\eta(0)=0$ and
\[
  G^{-1} (0) \cap W	 = \{ (\Lam, 0) \;  (\Lam, 0) \in W \}	
\cup \{  (\zeta(\a), \, \a \us + \a\eta(\a)) \; |\a| < a \} . 
										\tag{\eqgi}
\]
}
\vspace{- 2.5 mm}
\n
{\bf \remconc}. To show that  the Hopf bifurcation actually  
occurs for given concrete
\linebreak 
examples, {\thmbasic} 
 is often more  practical than {\propmain}.
 See e.g. [\refKsome].	\END
\propspace

				In what follows in this section, we use the same notation
 in Section 2.
We set
$Y_1 := \{ u \otimes c_1 + v \otimes s_1 \; u, v \in V \}$
   as a closed subspace of $Y$. 
   We define 
$L_1 : \Vc \to Y_1$ by
$L_1\psi := \Realpart (\psi \otimes e_1)$ 
for any $\psi \in \Vc$
 and 
$T_1 : X_1 \to Y_1$  by
$T_1 w := \dot{w} - Aw$ for any 
$w \in X_1$. 
    Namely, 
\[
L_1(a + i b) = a \otimes c_1 - b \otimes s_1
\moji{for any} a, b \in V,		\tag{\eqlon}
\]
\[
T_1 (a \otimes c_1 + b \otimes s_1)
= (b - A a)\otimes c_1 - (a + A b)\otimes s_1
		\moji{for any} a, b \in U. 		\tag{\eqtoa}
\]
 In view of (\eqlon) the following result clearly holds:
\medn
{\bf \propisom}. {\sl  
(i) The operator
 $L_1$ is  isomorphic as a real linear operator
from  $\Vc$ to  $Y_1$. Here, we regard $\Vc$ as a real linear space.

{\vskip 0.2cm}
\n
(ii) The operator
 $L_1 |_{\Uc} \: \Uc \to X_1$ is  isomorphic as a real linear operator
from  $\Uc$ to  $X_1$. Here, we regard $\Uc$ as a real linear space.
}

\propspace

\n
{\bf \propkernel}. {\sl  
(i) \, $L_1\, \kernel (i - \Ac) =  \kernel(T_1)$,
\medn
(ii) \, $L_1 \range (i - \Ac) =  \range(T_1)$.
}
\medbreak
			{\it Proof.} 
We easily verify from  {\propisom} (ii) that if $w \in X_1$ then
there exists a unique $\psi\in U_c$ such that $w=L_1 \psi$ and
 $T_1 w = L_1(i - \Ac)\psi$.

		(i) Let $w \in L_1 \kernel(i - \Ac)$.
Then, $w \in X_1$. So, there exists a unique $\psi\in U_c$ 
such that $w=L_1 \psi$, which leads to
 $T_1 w = L_1(i - \Ac)\psi$. Since $L_1$ is one to one,
$\psi \in \kernel(i - A_c)$. Therefore, 
$T_1 w = L_1 0 = 0$.  So, $w \in \kernel (T_1)$
and $L_1 \kernel(i - \Ac) \subset \kernel(T_1)$.

			Conversely, let $w \in \kernel (T_1)$.
Then,  there exists a unique $\psi\in U_c$ 
such that $w=L_1 \psi$.  It follows that $L_1(i - \Ac)\psi = T_1 w = 0$. 
Since $L_1$ is one to one,
$(i - \Ac)\psi = 0$. Therefore,
$\psi  \in \kernel(i - A_c)$ and  
 $w  \in L_1 \kernel (i - \Ac)$.
 We conclude that $\kernel(T_1) \subset L_1 \kernel(i - \Ac)$.
		
		(ii) Simple argument by linear algebra leads to the
desired conclusion, as in the proof of  (i). 
So, we leave the proof to the reader.		
\END
\eqspace
{\bf \proppsd}.  {\sl  Let $\psi \in \Uc$ and 
 $w= L_1\psi$.
 \medn
 (i) {\hskip 0.2cm} $L_1 (i \psi)= \dot{w}$,
 \medn
 (ii) If $\psi \in \kernel(i - \Ac)$,  then 
 $\,L_1 (i \psi)= A w$.
 }
  \medbreak
 			{\it Proof}. (i)  
 $L_1(i \psi)= \text{Re\,}\biggl[  \dfrac{d}{dt} 
 (\psi  \otimes e_1 )\biggr]
 = \dfrac{d}{dt} \,\text{Re\,} (\psi  \otimes e_1) = \dot{w}$.
\par
		(ii) We immediately obtain the desired conclusion from (i)
and {\propkernel} (i).	\END

\newchapter
 
 \centerline{\bf 4. Proof of {\thmmain}}
\medbreak 
		Let $X$ and $Y$ be  real Hilbert spaces defined by (\eqpersp).
We denote the $n$-th Fourier coefficient  of $\varphi \in Y_c$ by
\[
{\hat \varphi}(n):= \dfrac{1}{2\pi}\int_{0}^{2\pi} e^{-int} \varphi(t) dt  {\hskip 0.4cm} (n\in \Z). 			\tag{\eqfrac}
\]
We set 
\[
X_0 :=U 
\moji{and} X_{\infty} := 
\{ \varphi \in X \; {\hat \varphi}(n)=0 
\moji{for} n=-1, 0, 1 \}		\tag{\eqxn}
\]
  as closed  subspaces  of $X$,
\[
Y_0 :=V
\moji{and} Y_{\infty} := 
\{ \varphi \in Y \; {\hat \varphi}(n)=0 
\moji{for} n=-1, 0, 1 \}					\tag{\eqyn}
\]
as  closed subspaces of $Y$.
Let $X_1$ (resp. $Y_1$) be a closed subspace of $X$ (resp. $Y$)
defined in Section 2 (resp. Section 3).
\medbreak
			{\it Proof of {\thmmain}}. \,
			We apply {\thmbasic}.
We use the notation in Sections 2 and 3.
We denote $\Lam= (\lam, \sigma) \in K \times \R$.
We define $g \in C^2(K \times \R \times X, Y)$ by
$g(\Lam, u)= u_t - (\sigma +1) f(\lam, u)$,
where $f(\lam, u):= A u + h(\lam, u)$.
By the assumption (\Bone) in Section 2, there exists 
$\,\psi_{\star} \in \kernel(i - \Ac) - \{ 0 \}$.
Then, $\Realpart \psis$ and $\Imaginarypart \psis$ are linearly
independent in $U$. 
So, there exists $d \in U$ such that
$(d_c, \psis)_{\Uc} = 1$.
 We define 
 $l = (l^1, l^2) \in \linear(X, \R^2)$ by
\[
l^1 u  := \frac{1}{\pi}\int_0^{2\pi} (d, u(t))_U \cos t \,dt
\,\, \moji{and}	\,\,
		l^2 u  := \frac{1}{\pi}\int_0^{2\pi}  (d, u(t))_U  \sin t \,dt
\]
for $u \in X$.
We set
$u_\star := L_1 \psi_{\star}
= \Realpart(\psi_\star \otimes e_1) \in X_1$.	
Then, $l u_{\star} = (1,0) = {\bm e}_1$.
Let $H \: K\times \R\times X \to \R^2\times Y$ 
be the bounded operator defined by (\eqH).
Then, by  (\KtwoFourth) and {\propkernel} (i),
$H(0, u_{\star}) = 
(l u_{\star} - {\bm e}_1, \, T_1 \us) = (0,0)$.
We set $DH^{\star}:= DH(0, \us)$. Then, we have
\[
DH^\star
\begin{pmatrix}
\lam  \\ 
\sigma  \\ 
u  \\ 
\end{pmatrix}
=
\begin{pmatrix}
l^1 u \\ 
l^2 u \\ 
u_t - Au - \sigma A \us - \lam f_{\lam u}^0 \us  \\ 
\end{pmatrix} ,	\tag{\eqmapdh}
\]
where $f_{\lam u}^0:= f_{\lam u}(0,0)$.
We verify that 
$S:= DH^\star|_{\R^2 \directsum X_0 \directsum X_1} :
\R^2 \directsum X_0 \directsum X_1 \to 
\R^2 \directsum Y_0 \directsum Y_1$ 
and 
$T:= DH^\star|_{X_{\infty} } :
X_{\infty} \to Y_{\infty}$ 
are well - defined by {\remcomm}
and that $DH^\star = S \directsum T$.
We note that
$T u = u_t - A u$ for any $u \in X_{\infty}$.
In view of  the below {\lemsi}
and {\lemti}, $DH^{\star}$ is bijective.
So,  by {\thmbasic} 
 $(\lam, u)= (0, 0)$
is a Hopf bifurcation point and 
 there exist an open neighborhood $W$ of (0,0)  
 in $\R^2 \times X$,
$a  \in (0, \infty)$ and  functions 
$\zeta \in C^1((-a, a), \R^2)$, 
$\eta \in C^1((-a, a), Z)$ such that $\zeta(0)=0$, $\eta(0)=0$ and
(\eqgi) holds. 
Here, $Z:= \{ u \in X \; l^1 u=0 \}$. 
So, (a) holds. 

			Next, we show the following (4.5) in preparation
to prove (b) and (c).
\[
\zeta(- \a) = \zeta(\a) 
\moji{and}
 \eta(- \a) = - \tau_\pi (\eta(\a)) 
		\moji{for any} \a \in [0, a). 	\tag{\eqgam}
\]
We set $U(\a):= \a \us + \a \eta(\a) \in X$ for any 
$\a \in (-a, a)$.
We define $V(\a) \in X$ by $V(\a):= \tau_\pi (U(\a))$.
Let  $\gamma \in (0, a)$ be a constant such that 
$\{ (\zeta(\a), V(\a)) \; \a \in [0, \gamma) \} \subset W$.
Then, $(\zeta(\a), V(\a)) \in G^{-1}(0) \cap W$
for any $\a \in [0, \gamma)$. 
So,  by {\thmbasic} for any $\a \in [0, \gamma)$
 there exists $\b \in (-a, a)$ such that
  $(\zeta(\a), V(\a)) = (\zeta(\b), U(\b))$.
  On the other hand, $l^1 V(\a) = - \a$ and $l^1 U(\b) = \b$.
  Therefore, $(\zeta(- \a), U(- \a)) = (\zeta(\a), V(\a))$
  for any $\a \in [0, \gamma)$.
  Actually, we easily verify from commonly used argument by 
  contradiction that 
\[
a = \sup\,\{ q \in (0, a) \; (\zeta(- \a), U(- \a)) = (\zeta(\a), V(\a)) 
  	\moji{for any} \a \in [0,q) \},
\]
 which implies (\eqgam). 
  
            By  (\eqgam), $\zeta'(0)=0$. So, (b) holds.
 		Finally, we show (c).
Let $\e$ be a positive constant
such that if $(\lam, \sigma, w) \in \R^2 \times X$ 
satisfies $|\lam| < \e, |\sigma| < \e$
and $\| w \|_{X} < \e$ then $(\lam, \sigma, w) \in W$.
Now, let $(\lam, v)$ be a solution of (\eqpde) of period 
$2\pi(\sigma + 1)$, $|\lam| < \e$, $|\sigma| < \e$,
${\tilde v} \in X$ and
$\|  {\tilde v} \|_X < \e$,  where
${\tilde v}(t):= v((\sigma + 1)t )$ for $t \in \R$.
For simplicity, we set $(p, q):= l {\tilde v}=(l^1 {\tilde v}, l^2 {\tilde v})$.
First we consider the case: $q = 0$.
Then $(\lam, \sigma, {\tilde v}) \in W$ is a  solution of 
$G(\Lam, u):=(l^2 u, g(\Lam, u)) = 0$.
By {\thmbasic} there exists $\a \in (-a, a)$ such that 
$(\lam, \sigma) = \zeta(\a)$ 
and ${\tilde v} = \a \us + \a \eta(\a)$.
If $\a < 0$ then 
${\tilde v} = \tau_\pi \{ (- \a) \us + (- \a) \eta(- \a) \}$
in view of (\eqgam) and $\tau_\pi \us = - \us$.
		Next, we consider the case: $q \not = 0$.
There exists   $\theta \in (0, 2\pi)$ such that 
$e^{i \theta}=(p - i q)/\sqrt{p^2+ q^2}$.
Then, $l^2 \tau_\theta {\tilde v} = 0$ and 
$(\lam, \sigma, \tau_\theta {\tilde v}) \in W$ 
is a solution of $G(\Lam, u)=0$.
So, the present case is reduced to the case: $q = 0$.
Therefore, (c) holds.	\END           
 \eqspace
 
		In the above proof, we use the following two lemmas:
 \eqspace
{\bf \lemsi}. {\sl  The operator $S$ is bijective.}
\eqspace
{\bf \lemti}. {\sl  The operator $T$ is bijective.}
\eqspace 

		{\it Proof of {\lemsi}}.  
By  (\Bone), {\remcomm} and the implicit function theorem
(see e.g. [{\refCRbif}, Theorem A])
$f_u(\lam, 0)$ has an eigenvalue $\mu(\lam) \in \C$
and an eigenfunction $\psi(\lam) \in U_c$ 
corresponding to $\mu(\lam)$ 
for any $\lam$ in a small open interval $K_1$ 
such that  
$0 \in K_1 \subset K$,  $\mu(0)= i$, $\psi(0)= \psi_{\star}$,
$\mu(\cdot) \in C^2(K_1, \C)$ and $\psi(\cdot) \in C^2(K_1, U_c)$.
Differentiating 
$f_u(\lam, 0) \psi(\lam)=\mu(\lam)\psi(\lam)$
with respect to $\lam$, we have
\[
\mu'(0)\psi_{\star} + i \psi'(0)	
	= f^0_{\lam u}\psi_{\star} + \Ac \psi'(0).		\tag{\eqmuz}
\]
We set $p:= \text{Re}\, \mu'(0)$ ($\not = 0\,$  by (\Btwo)), 
$q= \text{Im}\, \mu'(0)$ and
$u_{\sharp}:= L_1 \psi'(0) \in X_1$.
 It follows from (\eqmuz) and 
{\proppsd} that
\[
 f_{\lam u}^0 \us = p \us + q A \us + T_1 u_{\sharp}.	 \tag{\eqfly}
\]
Let $u_0 \in X_0$,  $u_1 \in X_1$ and $u= u_0 + u_1$. 
In view of  (\eqmapdh) and (\eqfly), we have
\[
S	
\begin{pmatrix}
 (\lam, \sigma) \\ 
u_0    \\
u_1      \\
\end{pmatrix}		
=
\begin{pmatrix}
 l u_1 \\ 
 - A u_0   \\
T_1(u_1 - \lam u_{\sharp}) - \lam p \us 
  - ( \sigma + \lam q) A \us      \\
\end{pmatrix}
.		\tag{\eqsvec}
\]
By (\Bone),  we have
$\range(i - \Ac) \directsum \span \{ \psis \} =  \Vc$.
It follows from  {\propisom} (i), {\propkernel} (ii) and {\proppsd} 
that
\[
\range(T_1)	\directsum \span\{ \us, A\us \}  = Y_1.	\tag{\eqrto}
\]

		First, we show that $S$ is one to one.
Let $S(\lam, \sigma, u)= 0$.
It follows from (\Btwo),  (\eqsvec), (\eqrto) and $0 \in \rho(A)$
that $u_0=0$, $\lam= \sigma=0$,
\[
l u_1=0 \moji{and} T_1 u_1=0.			\tag{\eqluo}
\]
Let $\psi_1 := L_1^{-1} u_1 \in U_c$. 
Then by (\eqluo)  and {\propkernel} (i),
\[
\psi_1 \in \kernel(i - \Ac) \moji{and} (d, \psi_1)_{U_c} = 0.		
													\tag{\eqpin}
\]
It follows from (\eqpin), (\Bone) and $ (d, \psis)_{U_c}  = 1$  that 
$\psi_1 = 0$, which implies $u_1=0$. 
So, $S$ is one to one. 

			Next, we show that $S$ is onto.
Let $(a, b, y_0, y_1) \in \R^2 \directsum Y_0 \directsum Y_1$.
In view of  $\,0 \in \rho(A)$,
there exists $x_0 \in X_0$ such that $- A x_0 = y_0$.
By (\eqrto) there exist $w \in \range(T_1)$ and 
$(\gamma, \delta) \in \R^2$ such that 
\[
w + \gamma u_{\star} + \delta A u_{\star} = y_1.		\tag{\eqwg}
\]
We set $\lam_0:= - \gamma/p$ and 
$\sigma_0 := -\delta + \gamma q /p$.
There exists $v_1 \in X_1$ such that 
$T_1 (v_1 - \lam_0 u_{\sharp}) = w$.
Let $(\a, \b) := l v_1 \in \R^2$ and
$x_1:= v_1 + (a-\a)u_{\star}+(\b - b)Au_{\star}$.
By {\propkernel} (i) and {\proppsd} (ii),
we have 
$A u_{\star} = L_1 (i \psi_{\star}) \in \kernel(T_1)$.
So, $lAu_{\star}=(0, -1)$.
It follows from $l \us = \eone$, {\propkernel} (i),
(\eqsvec) and (\eqwg) that
$S(\lam_0, \sigma_0, x_0, x_1)=(a,b,y_0,y_1)$.
Therefore, $S$ is onto.		\END
\eqspace
\medbreak
			{\it Proof of {\lemti}}.  
Let $\cX := (X_\infty)_c$ and $\cY :=  (Y_\infty)_c$
(i.e. $\cX$ and $\cY$ be the complexification of 
$X_{\infty}$ and $Y_{\infty}$, respectively.)
It suffices to show that $T_c \: \cX \to \cY$ is bijective.
 Let $z \in \cY$. Then, 
$z = \sum_{|n| \ge 2} p_n \otimes e_n$ in  $\cY$,
where $\,p_n := {\hat z}(n)$.
It suffices to show that the following equation (\eqwt)
has a unique solution in $\cX$:
\[
u_t - A_c u = z		\tag{\eqwt}
\]
If a solution of (\eqwt) exists, we  obtain formally
 from the Fourier analysis that 
$\,u= \sum_{|n| \ge 2}   q_n \otimes e_n$, 
where $\,q_n := (i n - A_c)^{-1} p_n$.	
The proof is complete if we  show $u \in \cX$, i.e.
\[
\sum_{|n| \ge 2} (|n|^2 \| q_n \|^2_{V_c} 
		+ \|q_n\|^2_{U_c}) < \i.			\tag{\eqndn}
\]
It follows from (\Kone) that
\[
|n|\| q_n \|_{V_c} 
	\le |n|\, \| (in - A_c)^{-1} \|_{V_c \to V_c}
			\| p_n \|_{V_c} \le M \| p_n \|_{V_c},		
											\tag{\eqnn}
\]
\[
\|q_n\|_{U_c} = \| A_c(in - A_c)^{-1} p_n \|_{V_c}
	= \| in (in - A_c)^{-1}p_n - p_n \|_{V_c}
	\le (M + 1) \| p_n \|_{V_c}.		
										\tag{\eqdnu}
\]
By (\eqnn), (\eqdnu) and  Parseval's identity  we have (\eqndn). 		\END

\newchapter

\centerline{\bf 5. Examples}

\propspace

	   In this section we freely use the notation used in Section 4.	
\propspace

\n
{\bf Example 1}.
We consider the following  Cauchy problem:
\[
\begin{cases} 
& {\hskip - 0.3cm} u_{t}= u_{xx}  - v
	- \rho u + 
		u (\lam\,\kappa^2 - u^2 - v^2) 
		\moji{for} (x, t) \in \R \times [0, \infty), \hfill \\ 
&  {\hskip - 0.3cm} v_{t}= v_{xx}  + u
	- \rho v + 
		v  (\lam\,\kappa^2 - u^2 - v^2) 
		\moji{for} (x, t) \in \R \times [0, \infty), \hfill \\ 
\end{cases}
			\tag{\eqhc}
\]
Here, $\rho$ and $\kappa$
are functions on $\R$ defined by 
$\rho(x) := \left\{2 \tanh^2(x/2) - 1  \right\} / 4$
and $\kappa (x):= \sech(x/2)$.
 
		For the equation (\eqhc) the branch of periodic solutions 
 $\,(u, v)={  (u_{\lam}, v_{\lam})}$ $(\lam > 0)$
\linebreak
bifurcates  at $\lam=0$ from the branch of trivial solutions.
Here, ${ u_{\lam}(x, t) := 
  \sqrt{\lam} \,\kappa(x) \cos t }$ and 
${  v_{\lam}(x, t) :=}$ 
 ${\sqrt{\lam} \,\kappa(x) \sin t \,}$.
 
		We can  not  apply
   [\refCRhopf, Theorem 1.11]   to (\eqhc) since 
the linear operator in (\eqhc) does not have compact resolvents.
We show in  what follows that
we can apply our {\propmain} to (\eqhc) to verify 
the above Hopf bifurcation.
	
		 Let $V := L^2(\R) \times L^2(\R)$. 
We define $A \: V \to V$ and $A_k \: V \to V$ $(k= 0, \i)$ by 
the following:
\[
A
\begin{pmatrix}
\phi  \\ 
\psi  \\ 
\end{pmatrix}
:=
\begin{pmatrix}
\phi_{xx}  - \psi - \rho \phi  \\ 
\psi_{xx}  + \phi - \rho \psi  \\ 
\end{pmatrix}
\moji{for} (\phi, \psi) \in \domain(A):= H^2(\R) \times H^2(\R),
\]
\[
A_0
\begin{pmatrix}
\phi  \\ 
\psi  \\ 
\end{pmatrix}
:=
\begin{pmatrix}
\phi_{xx}  - \psi   \\ 
\psi_{xx}  + \phi   \\ 
\end{pmatrix}
\moji{for} (\phi, \psi) \in \domain(A_0):= H^2(\R) \times H^2(\R),
\]
\[
A_{\i}
\begin{pmatrix}
\phi  \\ 
\psi  \\ 
\end{pmatrix}
:=
\begin{pmatrix}
\phi_{xx}  - \psi - \phi/ 4  \\ 
\psi_{xx}  + \phi - \psi / 4  \\ 
\end{pmatrix}
\moji{for} (\phi, \psi) \in \domain(\Ai):= H^2(\R) \times H^2(\R),
\]
We set $U := \domain(A)= H^2(\R) \times H^2(\R)$.
By the Sobolev embedding theorem, 
$H^1(\R)$ is embedded in $L^6(\R)$.
So,  we can well define the map 
$h \: \R \times U \to V$ by
\[
h(\lam, \psibold) := 
\begin{pmatrix}
\phi (\lam\,\kappa ^2 - \phi^2 - \psi^2) \hfill \\
\psi (\lam\,\kappa ^2 - \phi^2 - \psi^2) \hfill \\
\end{pmatrix}
\moji{for} \lam\in \R \moji{and} 
\psibold=(\phi, \psi) \in U.		\tag{\eqhlp}
\]
Thus, (\Ktwo\ - 1) holds.

			We explain the correspondence between the present example
and the description in Section 4. We define $\usbold \in X_1$ by
$\usbold(x, t):= (\kappa(x)\cos t, \kappa(x)\sin t)$
and $\psisbold \in U_c$ by 
\[
\psisbold := L_1^{-1} (\usbold) = 
			(\kappa, - i \kappa).		\tag{\eqsp}
\]
We verify that $\usbold \in \kernel(T_1)$ and
$\psisbold \in \kernel(\Ac - i)$ (see  {\propkernel} (i)).
We set $\lam(\a)= \a^2$, $\sigma(\a) \equiv 0$ 
and $\etabold(\a) \equiv (0, 0)$ for $\alpha\in\R$.
Then,  we verify that
$(\lam, \sigma, u) = (\lam(\a), \sigma(\a), \a \usbold + \a \etabold(\a))$
is a solution of (\eqtpe). 
		
		In order to verify the applicability of our {\propmain} to (\eqhc) 
we describe the outline of the derivation of
 (\Kone), (\Ktwo), (\Bone)-(\Bthree)
 in what follows.

\propspace
\n
\underline{Derivation of   (\Ktwo)}
\medn	
We showed that (\Ktwo\ - 1) holds.
Let $\cV$ be a Hilbert space. 
We set $\cX:= H^1 ((0, 2\pi), L^2(\R))$
$\cap\, L^2((0, 2\pi), H^1(\R))$ and
 $\cY:= L^2((0, 2\pi), L^2(\R))$.
We use the next  embedding inequality
\linebreak
 (\equni) from
$H^1 ((0, 2\pi), \cV)$ into $C([0, 2 \pi], \cV)$ and
 the Gagliardo-Nirenberg  interpolation
\linebreak
 inequality (\eqint):
\begin{gather*}
 \| u \|_{C([0, 2 \pi], \cV)} 
		\le C_1 \| u \|_{H^1 ((0, 2\pi), \cV)}	 		
	\moji{for any} u \in  H^1 ((0, 2\pi), \cV),	 	\tag{\equni}  	\\
 \| \phi \|_{L^6(\R)} \le 
 \| \phi \|_{L^2(\R)}^{2/3}  \| \phi' \|_{L^2(\R)}^{1/3}
 		\moji{for any} \phi \in H^1 (\R).	 	\tag{\eqint}		\\
\end{gather*}
{\vskip - 0.57cm}	\noindent
Here, $C_1 >0$ is a certain constant independent of $u$.	
It follows from (\eqint), (\equni) and H{\"o}lder's inequality that 
if $f,g,h \in \cX$ then $f g h \in \cY$ with the estimate
\[
\|f g h  \|_{\cY}	 \le
C_2  \| f \|_{\cX} \| g \|_{\cX} \| h \|_{\cX}	
\moji{for any} f, g, h \in \cX.			\tag{\eqfgh}
\]
Here, $C_2 := C_1^{2/3} > 0$.
In view of (\eqfgh), we have
 (\Ktwo\ - 2).
In the present case the map $\Phi \: \R \times X \to Y$ is 
defined  by
\[
\Phi(\lam, {\bm u}) := 
\begin{pmatrix}
u (\lam\,\kappa ^2 - u^2 - v^2)   \hfill \\
v (\lam\,\kappa ^2 - u^2 - v^2)   \hfill \\
\end{pmatrix}
\moji{for} \lam\in \R \moji{and} 
{\bm u}=(u, v) \in X.				\tag{\eqplp}
\]
We omit the derivation of
 (\Ktwo\ - 3) since it is not difficult 
 by using (\eqfgh).	Finally,  (\Ktwo\ - 4) clearly holds. 
 	\END

{\vskip 0.3cm}
\propspace

		Next, we make preparation to derive (\Bthree) and (\Kone):
\n
\propspace
{\bf \lemfol}. 
{\sl  The following holds:
\propspace
\n
(i) $k=0$ or $|k| \ge 2 \,\Longrightarrow\,  i k \in \rho(\Azc)$
and $\| (\Azc - i k)^{-1} \|_{\Vc \to \Vc} \le 1$.
\medn
(ii)  $|k| \ge 4 \,\Longrightarrow \,
		\| (\Azc - i k)^{-1} \|_{\Vc \to \Vc} \le \dfrac{\sqrt{2}}{|k|}$.
}	
\propspace

			{\it Proof.} 
	Since $A_0$ is a real operator, it suffices to prove (i) and (ii)
	only for $k \ge 0$.
			
		(i)  Let $k=0$ or $k \ge 2$.
We consider the following equation:
\[
(\Azc - i k) (\phi, \psi) = (\gamma, \omega).		\tag{\eqazi}
\]
For any given $(\gamma, \omega) \in V_c\,$ this equation has a unique solution
$(\phi, \psi) \in U_c$ such that
\[
  {\hat \phi}(\xi) = \frac{- (\xi^2 + i k) {\hat \gamma}(\xi) 
	+ {\hat \omega}(\xi)}{1 + (\xi^2 + i k)^2}
\moji{and} 	
{\hat \psi} (\xi)= - \frac{ {\hat \gamma}(\xi) 
	+ (\xi^2 + i k) {\hat \omega}(\xi)}{1 + (\xi^2 + i k)^2}.
						\tag{\equh}	
\]
By  Schwarz inequality
\[
|{\hat \phi}(\xi)|^2 + |{\hat \psi}(\xi)|^2 
= \frac{ (\xi^4 + k^2 + 1)(|{\hat \gamma}(\xi)|^2 + |{\hat \omega}(\xi)|^2)
	+ 4k \Realpart \{ i \,\overline{ {\hat \gamma}(\xi)} \,{\hat \omega}(\xi)  \} }
				{\xi^8 + 2(k^2 + 1)\xi^4 + (1 - k^2)^2}  		\tag{\equhz}
\]	
\[
{\hskip - 0.5 cm}
\le 
	J(k, \xi)	\{ |{\hat \gamma}(\xi)|^2 + |{\hat \omega}(\xi)|^2 \},	
	\moji{where} J(k, \xi) :=  
\frac{\xi^4 + (k + 1)^2}{\xi^8 + 2(k^2 + 1)\xi^4 + (1 - k^2)^2}.
\]
We verify that $J(k, \xi) \le 1$ for any $\xi\in\R$.
So, we have the desired result.		
\propspace		

		(ii)  Let $k \ge 4$. Then, we verify that
 $J(k, \xi) \le 2/k^2$ for any $\xi \in \R$.
By this and (\equhz), we have the desired result.		\END

\propspace
\n
\underline{Derivation of (\Bthree) and  (\Kone)}
	
\n
We define $B \in \linear(V)$  by 
$B(\phi, \psi) = 
      - (\rho \phi, \rho \psi)$ for   $(\phi, \psi) \in V$.
Then,  we have
\linebreak
$B_c =  (A_c - i k) -  (\Azc - i k)$ on $\domain(\Azc) \,
(= \domain(A_c)\,)$.
By {\lemfol} (i) we have
$i k \in \rho(\Azc)$ and
$\|B_c  \|  \| (\Azc - i k)^{-1} \| \le 1/4 \, ( < 1)$ for any
$k \in \Z - \{ \pm 1 \}$.	
So, (\Bthree) holds in view of stability property of 
bounded inverse operator
(see e.g. [\refKcom, Corollary 2.4.1]).     
 		
		Next, by {\lemfol} (ii) we have
$\| (\Azc - i k)^{-1} \|_{\Vc \to \Vc} \le 1/2 \sqrt{2}  < 1/2$
for $k \ge 4$.
It follows from [\Ktwo, Corollary 2.4.1] that
\[
\| (A_c - i k)^{-1} \| 
\le \frac{\| (\Azc - i k)^{-1} \|}{1 - \| B_c \|  \| (\Azc - i k)^{-1} \|}
\le \frac{8 \sqrt{2}}{7 k} \moji{for} k \ge 4.
\]
In view of this estimate and (\Bthree), we obtain (\Kone).		\END

{\vskip 0.3cm}
\propspace

		We make preparation to derive (\Bone).

\propspace
\n
{\bf \lembas}. 
{\sl  
(i)  $\Ac - i$ is a Fredholm operator  of index 0.
\medn
(ii) $\kernel(\Ac - i) = \span\{ \psisbold \}$.
}
\propspace
           
			{\it Proof.}  (i) 
Step 1.  By using the Fourier analysis as in the proof of {\lemfol} (i),
we verify that for any $(\gamma, \omega) \in V_c$
the equation $(\Aic - i) (\phi, \psi) = (\gamma, \omega)$
 has a unique solution $ (\phi, \psi) \in U_c$.
So, $\Aic - i$ is bijective. 

		Step 2. By Step 1, $\Aic - i$ is a Fredholm operator of index 0.
We define ${\hat A}, {\hat A}_{\i} \in \linear(U, V)$ by 
${\hat A}{\bm \phi} = A {\bm \phi}$ for ${\bm \phi} \in U$ and by
${\hat A}_{\i} {\bm \phi} = \Ai {\bm \phi}$ for ${\bm \phi} \in U$.
 In order to complete the proof
 it suffices to show that ${\hat A}_c - i$ is a Fredholm operator 
of index $0$. 
We denote by $\chi_R$  the identity function of 
$(- R, R)$ for $R > 0$. 
Let $g(x) := -  \{ 1 - \tanh^2 (x/2) \} / 2$.
We define $H_R \in \linear(U_c, V_c)$ by
 $H_R := ({\hat A}_{\i c} - i) - (1 - \chi_R)g$ for $R > 0$.
Then, ${\hat A}_c - i = H_R - \chi_R g$.
In view of stability property of bijective operator
 (e.g. [\refKcom, Corollary 2.4.1])
$H_R$ is bijective for a sufficiently large constant $R > 0$.
Since the  multiplication operator 
$\chi_R g \,\cdot : U_c \to V_c$ is compact,
${\hat A}_c - i$ is a Fredholm operator of index 0
by the stability property of Fredholm operator.

			(ii)  We consider the following equation:
\[
(A_c - i) (\phi, \psi) = (0, 0).		\tag{\eqazz}
\]
Let ${\cal H} := 
\{ (\phi, \psi) \in \{ C^2(\R) \}^2 \; (\eqazz) \text{ holds} \}$
be a complex linear space.
We note that $\kernel(A_c - i) \subset {\cal H}$ and 
$\kernel(A_c - i) = {\cal H} \cap U_c$.
Let ${\cal H}_1$ and ${\cal H}_2$ be subspaces of ${\cal H}$
defined by 
\[
{\hskip -1.5cm}
  {\cal H}_1 := \{  (\phi, - i \phi) \in {\cal H} \; 
\phi'' - \rho \phi = 0 \moji{on} \R \},	
{\hskip 0.2cm}
 {\cal H}_2  	:= \{  (\phi, i \phi) \in {\cal H} \; 
\phi'' - (2i + \rho) \phi = 0  \moji{on} \R  \}.	
\]
Then we verify $\dim {\cal H} = 4$ and
 $\dim {\cal H}_1 = \dim {\cal H}_2 = 2$
 from the foundational theorem on
 \linebreak 
 existence and uniqueness of solutions for ODEs
 and standard techniques on ODEs.
 \linebreak
 Actually we have ${\cal H} = {\cal H}_1 \directsum {\cal H}_2$
 since  four vectors $(1, 0, \pm i, 0)$ and $(0, 1, 0, \pm i)$
are linearly independent and two solutions satisfying
$(\phi(0), \phi'(0), \psi(0), \psi'(0)) = (1, 0, - i, 0), (0, 1, 0, - i)$
 (resp. $(1, 0,  i, 0), (0, 1, 0, i)$) belong to ${\cal H}_1$ 
  (resp. ${\cal H}_2$).
We note that ${\bm \phi} = (\phi, -i \phi) \in {\cal H}_1$,  
${\bm \psi} = (\psi, i \psi) \in {\cal H}_2$ and 
${\bm \phi} + {\bm \psi} = (\phi + \psi, i(\psi - \phi)) \in U_c$
implies $\phi, \psi \in H^2(\R)$ and ${\bm \phi}, {\bm \psi} \in U_c$.
Combining this and $\psisbold \in {\cal H}_1 \cap U_c$,
the proof is complete if we show that ${\cal H}_1 \not \subset U_c$
and ${\cal H}_2 \cap U_c = \{ 0 \}$. 
First, we show  ${\cal H}_1 \not \subset U_c$.
To this end, it suffices to show that 
there exists ${ \bm z} =(z, -i z)\in {\cal H}_1$ such that 
\[
z(x)  \ge w(x) := e^{x/2} +1	 \moji{for any} x \ge 6.		\tag{\eqzin}
\]
We  verify that
$w'' (x)- \rho (x) w (x)  
=  - (1/16)(e^{x} + e^{- x} - 8 e^{x/2} - 6) \,\sech^2 (x/2) < 0
$ for any $x \ge 6$. 
Let $z$ be a solution of $z'' - \rho z =  0$ 
satisfying the initial condition:
$z(6) = w(6)$ and $z'(6) > w'(6) = e^3 / 2$.
Then, $\phi := z - w$ satisfies 
$\phi'' - \rho \phi \ge 0 \moji{on} [6, \i).	$
Let $a \in (6, \i)$. By $\phi'(6) > 0$ and the maximum principle, 
$\phi(y)$ achieves the maximum value at $y = a$
 on the interval $[6, a]$. 
So, $\phi(y)$ is actually monotone increasing for $y \ge 6$ and
(\eqzin) holds. Therefore, 
${\cal H}_1 \not \subset U_c$.

				Next, we show ${\cal H}_2 \cap U_c = \{ 0 \}$.
We define the operators 
$H$, $H_{\i} \: L^2(\R) \to L^2(\R)$ by
\[
 H \phi := \phi'' - (2i + \rho) \phi,	{\hskip 0.4cm}
 H_{\i} \phi := \phi'' - (2i + 1/4) \phi	
\]
for $\phi \in \domain(H)= \domain(H_{\i}) := H^2(\R)$.
We verify
from the standard Fourier analysis  that 
for any given $\psi \in L^2(\R)$  the equation $H_{\i} \phi = \psi$ 
has a unique solution
$\phi \in H^2(\R)$ such that 
$ {\hat \phi}(\xi) = - \, {\hat \psi}(\xi) / (2 i + 1/4 + \xi^2)\,$
for any  $\xi \in \R$.
So, $H_{\i}$ is bijective, $\| H_{\i}^{-1} \| \le 1/ \sqrt{(1/4)^2  + 4}$
and $\| B \| \le 1/2$. Here,
$B \in \linear(L^2(\R))$ is the multiplication operator
satisfying $B = H - H_\infty$ on $\domain(H)$.
In view of $\| H_{\i}^{-1} \| \| B \| < 1$ and
and the stability of inverse operators
 (e.g. [\refKcom, Corollary 2.4.1]) $H$ is also bijective.
  Therefore, ${\cal H}_2 \cap U_c = \{ 0 \}$.		\END
 
\pagebreak
\noindent
\underline{Derivation of (\Bone)}
\medn
It suffices to show that $i$ is the simple eigenvalue of $\Ac$. 
In view of {\lembas}, 
the proof is complete if we show $\psisbold \not \in \range(\Ac - i)$.
We proceed by contradiction. 
Suppose $\psisbold  \in \range(\Ac - i)$.
We verify that $\psisbold \in \kernel(\Ac^* + i)$, which implies 
$(\psisbold, \psisbold)_{\Vc} = 0$.
This contradicts $\psisbold  \in \Vc - \{ 0 \}$. 		\END	

\propspace
\n
\underline{Derivation of (\Btwo)}
\medn	
As in the proof of  {\lemsi}, it follows from (\Bone), {\remcomm}
and the implicit function theorem that
  $A_c + \lam \kappa^2$ has an eigenvalue 
  $\mu(\lam) \in \C$
and eigenfunction $\psibold(\lam) \in \domain(A_c)$
corresponding to $\mu(\lam)$  
for any $\lam$ in a small  neighborhood of $0$ 
such that   
$\mu(\cdot)$ and $\psi(\cdot)$
 are functions  of class $C^2$
 with $\mu(0)= i$ and $\psibold(0)= \psisbold$.
It follows from (\eqmuz) that
$
\mu'(0) \psisbold =	
	f^0_{\lam u}\psisbold + (\Ac - i) {\bm \psi}'(0).		
$
Combining this and $\psisbold \in \kernel ( (\Ac - i)^*)$,
\[
\mu'(0) \| \psisbold \|_{V_c}^2 = (\psisbold, f^0_{\lam u} \psisbold)
= 2 \int_{\R} \sech^4 (x/2) \, dx >0.
\]
So, $\mu'(0) > 0$.		\END

\propspace
{\vskip 0.5cm}

\n
{\bf Example 2}.
We consider the following system of Fitzhugh-Nagumo type:
\[
\begin{cases}
& u_{t}= v - u    
		\moji{for} (x, t) \in I \times [0, \infty), \hfill \\ 
& v_{t}= v_{xx} - 2u + 2v 
	+ u (\lambda \sin^2 x - 2u^2 + 2 uv- v^2) 
		\moji{for} (x, t) \in I \times [0, \infty), \hfill \\ 
&  u(0, t)= u(\pi, t)=0 \moji{for} t \in  [0, \infty), \hfill	\\ 
&  v(0, t)= v(\pi, t)=0 \moji{for} t \in  [0, \infty).  \hfill
\end{cases}
			\tag{\eqfns}
\]
Here, we set $I := (0, \pi)$.
The branch of periodic solutions 
 $\,(u, v)={  (u_{\lam}, v_{\lam})}$ $(\lam > 0)$
bifurcates  at $\lam=0$ from the branch of trivial solutions.
Here, $u_{\lam}(x, t) := \sqrt{\lam} \cos t \sin x$ and 
$v_{\lam}(x, t) := \sqrt{\lam}(\cos t - \sin t)\sin x$.
Let $V := H^1_0 (I) \times L^2(I)$. We define $A \: V \to V$ by 
  \[
{\hskip - 0.7 cm}
A \ubold := 
 \begin{pmatrix}
v - u  \\ 
v_{xx} -2 u + 2 v  \\ 
\end{pmatrix}
\moji{for any}
\ubold :=
\begin{pmatrix}
u  \\ 
v  \\ 
\end{pmatrix}
 \in 
	U= \domain(A):= H^1_0(I) \times (H^1_0 \cap H^2)(I).
\]
For any $z \in \rho(A)$,
the resolvent $(z - A)^{-1}$ is not compact.
Actually, we have
\[
(z - A)^{-1} 
	\biggl (\frac{z + 1}{n} \sin n x, \,\, - \frac{2}{n}  \sin n x \biggr )
	= \biggl ( \frac{\sin n x}{n}, \,\, 0 \biggr )
	\moji{for any}  n \in \N.
\]
So, $(z - A)^{-1}$ maps a bounded sequence in V to a
non-precompact sequence.
Therefore, we can not apply [\refCRhopf, Theorem 1.11] to 
our problem (\eqfns).

We show that we can apply our {\propmain} to verify 
the above Hopf bifurcation.
To this end, we describe the outline of the derivation of
 (\Kone), (\Ktwo), (\Bone)-(\Bthree) in what follows.
In the same way as Example 1, we can well define the map
$h \: \R \times U \to V$ by
\[
h(\lam, \ubold)
:=(0, u (\lambda \sin^2 x - 2u^2 + 2 uv- v^2))
\moji{for any} \lam \in \R \moji{and} \ubold := (u, v) \in U.
\]
Thus, (\KtwoFirst) holds.

			Let  $X$ and $Y$ be  real Hilbert spaces defined by (\eqpersp).
We define $\usbold \in X_1$ by
$\usbold(x, t):= (\cos t \sin x, (\cos t - \sin t)\sin x)$
and $\psisbold \in U_c$ by 
$\psisbold(x, t):= L_1^{-1}(\usbold)= (\sin x, (1 + i)\sin x)$.
Then, we verify that $\usbold \in \kernel(T_1)$ and
$\psisbold \in \kernel(\Ac - i)$ (see  {\propkernel} (i))
and that 
$(\lam, \sigma, u) = (\lam(\a), \sigma(\a), \a \usbold + \a \etabold(\a))$
is a solution of (\eqtpe). 
Here, we set $\lam(\a)= \a^2$, $\sigma(\a) \equiv 0$ 
and $\etabold(\a) \equiv 0$ for $\alpha\in\R$.
\propspace
\n
\underline{Derivation of   (\Ktwo)}  
		
\medn	
We showed that (\Ktwo\ - 1) holds.
We set $\cX:= H^1 ((0, 2\pi), L^2(I)) \cap 
L^2((0, 2\pi), H^1_0(I))$ and $\cY:= L^2((0, 2\pi), L^2(I))$.
It follows from (\eqint) that
\[
\| \phi \|_{L^6(I)} \le 
 \| \phi \|_{L^2(I)}^{2/3}  \| \phi' \|_{L^2(I)}^{1/3}
 		\moji{for any} \phi \in H^1_0 (I).	 	\tag{\eqini}
\]
In view of (\equni), (\eqini) and H{\"o}lder's inequality,
if $f,g,h \in \cX$ then $f g h \in \cY$ with the estimate (\eqfgh).
By using (\eqfgh), we verify that (\Ktwo\ - 2) and (\Ktwo\ - 3)
hold. Clearly, \linebreak
(\Ktwo\ - 4) holds.
		\END

\propspace
\n
\underline{Derivation of   (\Bone)}  
	\medn	
Let $(u, v) \in U_c$ satisfy $(i - A_c)(u, v) = (0, 0)$.
Then,  eliminating $v$ from this equation, we have 
$u_{xx} + u = 0$. It follows that $\kernel (i - A_c) = \span \{ \psisbold \}$. Now, it suffices to show the following:  
\propspace
(\eqiac)  {\hskip 0.2cm}  For any ${\bm p} \in V$ there exist 
$\ubold \in U_c$ and a unique $\alpha \in \C$ such that
{\vskip - 0.9 cm}
\[
 {\hskip 1.5 cm}	(i - A_c) \ubold + \alpha \psisbold = {\bm p}.
\]
Let $\ubold := (u, v) = (\sum u_n \sin nx, \sum v_n \sin nx)$ and
 ${\bm p} := (p, q) =  (\sum p_n \sin nx, \sum q_n \sin nx)$
 be Fourier expansion of $\ubold$ and ${\bm p}$. 
 Here, $u_n$, $v_n$, $p_n$ and $q_n$ are complex numbers.
 Then, by elementary Fourier analysis, 
 the solutions of the equation in (\eqiac) are given by the following:
 \[
u= \dfrac{v + \alpha \sin x - p}{1 + i}, 
	\q v= c \sin x + \sum_{n= 2}^\infty 
		\dfrac{2p_n - (1 + i) q_n}{(i + 1)(n^2 - 1)} \sin n x	
							 \tag{\equva}
\]
{\vskip - 1.0 cm}
\[
{\hskip 9.0 cm} \moji{and} 
  		\alpha= \dfrac{(1 + i) p_1 - i q_1}{2},	
\]
where $c$ is any complex number. So, (\eqiac) holds.  \END

\propspace
\n
\underline{Derivation of   (\Btwo)}  
\medn	
By the integration by parts, we verify 
${\bm \psi_{\#}} := ( (1+ i) \sin x, \sin x) \in \kernel ((i - A_c)^*)$.
In the same way as Example 1, 
$
\mu'(0) ({\bm \psi_{\#}}, \psisbold)_{V_c} 
		= ({\bm \psi_{\#}}, f^0_{\lam u} \psisbold)_{V_c}.
$
It follows that 
$\mu'(0) = (1/\pi) \int_0^\pi \sin^4 x \, dx  > 0$.	\END

\propspace
\n
\underline{Derivation of   (\Bthree) and (\Kone)}  
\medn	
It is similar to the derivation of (\Bone).
Let $k \in \N \cup \{ 0 \}$ and ${\bm a}= (a, b) \in V$.
We consider the equations $(A_c - k i) \ubold = {\bm a}$,
where $\ubold = (u, v) \in U$.
We set $\ubold := (u, v) = (\sum u_n \sin nx, \sum v_n \sin nx)$ and
$b - 2a/(1+k i)  = \sum d_n \sin nx$.
 Then, it follows from
\linebreak  
elementary Fourier analysis that
 \[
u= \dfrac{v - a}{1+ k i}	 \moji{and}
\left\{ - n^2 + \left( 2 - k i - \dfrac{2}{1 + k i}  \right)
\right \} v_n = d_n.		\tag{\eqvmp}
\]

		First, we consider the case: $k = 0$. 
In this case, $\ubold \in U$ is uniquely given by 
the first equality of (\eqvmp) and 
$v = - \sum(d_n/ n^2) \sin n x$. So, $0 \in \rho(\Ac)$.

		Next, we consider the case: $k \ge 2$. 
By (\eqvmp), $\ubold \in U$ is uniquely given by 
the first equality of (\eqvmp) and
\[
v_n = - \,\dfrac{d_n}{ \{n^2 - 2 + 2/(1+ k^2) \}  
		+  i k  \{ 1 - 2/(1+ k^2) \}  }.		\tag{\eqvne}
\]
It follows from (\eqvne) that $v \in H^1_0 \cap H^2(I)$, 
which leads to $\ubold \in U$. So, $k i \in \rho(\Ac)$ for
$k \ge 2$. Therefore, (\Bthree) holds.

		Finally, let $k \ge 2$ and we denote by $\| \cdot \|$
the norm of $L^2(I)$ and by $\| \cdot \|_1$ the norm of $H^1_0(I)$.  Here, $\| h \|_1 := \| h_x \|$ for $h \in H^1_0(I)$.
We easily verify $n |v_n| \le |d_n|$ for $n \in \N$.
So, $\| v_x \| \le  2\| a \|/ \sqrt{1+ k^2}
+ \| b \|$. By this and (\eqvmp), we have
\[
\| u_x \| \le \dfrac{1}{\sqrt{1+ k^2}}
	\left( \dfrac{2\| a \|}
		{\sqrt{5}} + \| b \| + \| a_x \| \right).
								\tag{\eqnux}
\]
On the other hand, by (\eqvne)
\[
|v_n| \le \dfrac{|d_n|}{k \{ 1- 2/(1+ k^2) \}} 
		\le \dfrac{5}{3k} |d_n| \moji{for} n \in \N.
\]
It follows that
\[
\| v \| \le \dfrac{5}{3k}
	\left( \dfrac{2\| a \|}
		{\sqrt{1+ k^2}} + \| b \|  \right).		
										\tag{\eqnov}
\]
Therefore, (\Kone) holds in view of (\eqnux), (\eqnov) and
Poincar{\'e} inequality 
$\| a \| \le \| a_x \|$.  	\END

\newchapter

\centerline{\bf Acknowledgments}
\medn
The author is grateful to Professor Masato Iida and
Professor Yoshitaka Yamamoto for useful discussion.

\newchapter

\centerline{\bf References}
\medn
[\refAma] H. Amann, Hopf bifurcation in quasilinear reaction--diffusion systems,
Delay \linebreak
Differential Equations and Dynamical Systems, Lecture Notes in Mathematics 1475 (1991) 53--63.
\medn
[\refBKST] T. Brand, M. Kunze, G. Schneider, T. Seelbach, Hopf bifurcation and exchange of stability in diffusive media, Arch. Rat. Mech. Anal. 171 (2004) 263--296.

\medn
[\refCRbif] M.G. Crandall and P.H. Rabinowitz,  Bifurcation from  
simple eigenvalues, J. Func. Anal. 8 (1971) 321--340. 
\medn
[\refCRhopf] M.G. Crandall and P.H. Rabinowitz, 
The Hopf bifurcation theorem in infinite
 \linebreak
  dimensions,
 Arch. Rat. Mech. Anal. {67} (1977) 53--72.
\medn
[\refGMW] D. Gomez, L. Mei, J. Wei,
Stable and unstable periodic spiky solutions for the Gray-Scott system and the Schnakenberg system, J. Dynam. Differential Equations 32 (2020), 441--481. 
\medn
[\refKsym] T. Kawanago, A symmetry-breaking bifurcation 
theorem and some related theorems
 applicable to maps having unbounded derivatives,
Japan J. Indust. Appl. Math. 21 (2004) 57--74.  \,\,
Corrigendum to this paper:  
Japan J. Indust. Appl. Math. 22 (2005) 147.
\medn
[\refKcom] T. Kawanago,  Computer assisted proof to 
symmetry-breaking 
bifurcation phenomena in nonlinear vibration,
Japan J. Indust. Appl. Math. 21 (2004) 75--108. 
  
 \medn
[\refKsome] T. Kawanago, Codimension-$m$ bifurcation theorems 
applicable to the numerical\linebreak 
verification methods, Advances in Numerical Analysis, 
vol. 2013 (2013), Article ID 420897.
\medn
[\refKie]  H. Kielh{\"o}fer, Bifurcation theory. An introduction with 
applications to partial
 \linebreak
  differential equations. Second edition. Applied 
Mathematical Sciences, 156. Springer, New York, 2012. 
\medn
[\refLiZY] H. Li,  X. Zhao and W. Yan, Bifurcation of time-periodic 
solutions for the \linebreak
 incompressible flow of nematic liquid crystals in three dimension, Adv. Nonlinear Anal. 9 (2020) 1315--1332.
\medn
[\refliu] Z. Liu, P. Magal and S. Ruan, Hopf bifurcation for 
non-densely defined Cauchy problems,
Z. Angew. Math. Phys. 62 (2011) 191--222.
\medn
[\refMS] A. Melcher and G. Schneider,
A Hopf-bifurcation theorem for the vorticity formulation of the Navier-Stokes equations in
$\R^3$,  Comm. Partial Differential Equations 33 (2008), no. 4-6, 772--783. 
\medn
[\refNishi] T. Nishida, Y. Teramoto and H. Yoshihara,
Hopf bifurcation in viscous incompressible flow down an inclined plane,
 J. math. fluid mech. 7 (2005) 29--71.
\medn
[\refWYZ] Q. Wang, J. Yang, L. Zhang, Time-periodic and stable patterns of 
a two-competing-species Keller-Segel chemotaxis model: effect of cellular growth, 
Discrete Contin. Dyn. Syst. Ser. B 22 (2017), 3547--3574.

\end{document}